\documentclass[12pt,a4paper]{article}

\usepackage{pdflscape}
\usepackage{amsmath}
\usepackage{amssymb}
\usepackage{latexsym}
\usepackage{srcltx}
\usepackage{graphics}
\usepackage{color}
\usepackage{epsfig}
\usepackage{color}
\usepackage{graphicx,epstopdf}

\newtheorem{thm}{Theorem}
\newtheorem{prop}{Proposition}

\newtheorem{cor}{Corollary}
\newtheorem{rem}{Remark}

\newtheorem{defn}{Definition}

\date{}

\begin{document}

\author{Vladimir Yu. Protasov 
	\thanks{Moscow State University;  {e-mail: \tt\small
			v-protassov@yandex.ru}
	}\ ,  
	Rinat Kamalov
	\thanks{Moscow Institute of Physics and Technology,   
		{e-mail: \tt\small rinat020398god@yandex.ruDISIM}
}}

\title{How do the lengths of switching   intervals\\ influence the stability of a dynamical system? 
}

\maketitle

\begin{abstract}
	
	If a linear switching system with frequent switches 
	is stable, will it be stable under arbitrary switches? 
	In general, the answer is negative. Nevertheless, this question can be answered 
	in an explicit form 
	for any concrete system. This is done  by finding  the mode-dependent critical lengths of switching intervals after which any enlargement does not influence the stability. 
	The solution is given in terms of  the exponential polynomials 
	of least deviation from zero on a segment (``Chebyshev-like'' polynomials).
	By proving several theoretical results on exponential polynomial approximation
	we   derive an algorithm for finding such polynomials  and
	for computing the critical switching time. The convergence of the algorithm is estimated and numerical results are provided.

	\bigskip

	\smallskip

	\begin{flushright}

	\end{flushright}
	
\end{abstract}
\bigskip

\vspace{1cm}

\begin{center}
	
	\large{\textbf{1. Introduction}}	
\end{center}
\bigskip 
\section{Introduction}
We consider linear switching systems of the form   $ \dot{\boldsymbol{x}}(t)=A(t)\boldsymbol{x}(t) $
with a given control set~$ \mathcal{A} $ and with mode-dependent dwell time restrictions~$ m(A) > 0, \, 
A \in {\mathcal{A}} $. 
Here the vector-function~$ \boldsymbol{x}(\, \cdot \,) $ is a trajectory, $ {\mathcal{A}} $ is a compact set of~$ d\times d $-matrices,~$ A(\, \cdot \,): \, [0, +\infty) \to \mathcal{A} $ 
is a piecewise-constant function called {\em switching law}. Each matrix from~$ \mathcal{A} $
is a {\em regime} of the system;  {\em switching points}
are points of discontinuity of~$ A(\, \cdot \,) $; the interval between two 
consecutive switching points is a {\em switching interval}. 
The  dwell time assumption 
means that every switching interval of the  regime~$ A\in \mathcal{A} $ is not shorter than
the prescribed length~$ {m(A) > 0} $.    

Such  dynamical systems with matrix control naturally appear  
in problems of robotics, electronic engineering, mechanics, planning, etc., see~\cite{GC}, \cite{GLP}, \cite{L}, \cite{MP} and references therein. 
A system is {\em asymptotically stable} if for every switching law and for every~$ \boldsymbol{x}_0 $, we have $\boldsymbol{x}(t) \to 0 $ as~$t\to +\infty $. 
Note that for linear switching systems, this property coincides 
with the locally attractiveness (see, for instance,~\cite[Theorem 3.4]{L}). 
We do not consider other types of stability and usually drop the word ``asymptotically''.

For one-element control sets~$ \mathcal{A} = \{A\} $, the stability is equivalent to that  the matrix  $ A $ is Hurwitz, i.e., all its 
eigenvalues have  negative real parts. If $ \mathcal{A} $ contains more than one matrix, then the stability problem becomes much 
harder~\cite{GC}, \cite{GO}, \cite{L}, although some practically applicable 
criteria of stability have been elaborated in the 
literature, see, for instance~\cite{BCM}, \cite{BM}, \cite{chesi0}, \cite{GLP}, \cite{SNM}. 
As a rule, they work for low dimensions, 
for higher dimensions they become rough (sufficient conditions are far from necessary ones). 
The stability depending on switching time intervals is a subject of an extensive 
research~\cite{chesi0}, \cite{CGPS}, \cite{ISW}, \cite{PK}, \cite{SFS}. 
We address the following  problem: {\em when the stability under frequent switches 
implies the stability without any  switching time restrictions}? In this case the stability can  be decided without applying the aforementioned criteria. 
Assume a linear system is stable under 
the assumption that the lengths of all switching intervals 
do not exceed certain mode-dependent bounds~$ M(A) $. 
This means that for every~$ \boldsymbol{x}_0$ and 
for every switching law~$A(\cdot)$ such  that the lengths 
of switching intervals of each mode~$A\in \mathcal{A}$ do not exceed~$M(A)$, we have $\boldsymbol{x}(t) \to 0 $ as~$t\to +\infty $. 
Will the system  remain stable after 
removal of those restrictions on the lengths? 
In practice, the frequent switching stability can either be established by sufficient conditions 
or be supposed empirically from experiments. The problem is to 
decide the stability of the original (unrestricted) system only from its frequent switching stability, without applying expensive and hard criteria.  

We solve this problem by finding the critical upper bounds for~$ M(A) $. 
To each Hurwitz  $d\times d$ matrix~$ A $, we  associate
a positive number~$ T_{cut}(A) $, which is   
the moment of time when a trajectory of the one-regime 
system~$ \dot{\boldsymbol{x}} (t)\, = \, A\, \boldsymbol{x}(t) $ enters its (symmetrized) convex hull. Theorem~\ref{th.10} proved in Section~3 states that if the stability takes place 
under the upper switching time  restrictions~$ M(A) = m(A) + T_{cut}(A) $,  i.e., if the switching time (after deducing the dwell time) is less than~$T_{cut}(A)$, then  
the system  stays stable without those restrictions. 

The dwell time assumption cannot be omitted: if for at least one regime~$ \bar A \in \mathcal{A} $, we have $ m(\bar A) = 0 $, then one can split all long switching intervals by momentary interfering with~$\bar A$, which  makes the presence of long switching intervals meaningless.   

Our main results establish a relation between the critical length~$T_{cut}$ and 
the exponential polynomials of lest deviations from zero on a segment (we call them 
``Chebyshev-like polynomials'').  Then we involve the tools of approximation theory
to find~$T_{cut}$.  This is done by 
extending the concepts of resonance and of the Remez algorithm to exponential polynomials o, more generally, to non-Chebyshev systems. This allows us to characterise the polynomial of lest deviations and 
then to find~$T_{cut}$ algorithmically.

The computation of the critical length~$ T_{cut} $ will be done in two steps. First we 
characterize this value geometrically, in terms of the convex hull of trajectory. 
Then we reduce its computation  to the problem 
of  finding the ``Chebyshev-like'' exponential polynomial.

To  formulate the main results we introduce the concept of {\em cut tail point}
of a matrix, which is probably of some independent interest. 

\section{Cut tail points}

For a Hurwitz matrix~$ A $, we consider a one-regime 
system~$ \ \dot{\boldsymbol{x}} \, = \, A\, \boldsymbol{x}, \ \boldsymbol{x}(0) = \boldsymbol{x}_0$. We denote by~$ L = L_{\boldsymbol{x}_0} $ the minimal by inclusion invariant subspace of~$ A $ that contains~$ \boldsymbol{x}_0 $. 
For a segment~$ [t_1, t_2] \subset \mathbb R_+$, let~$ \Gamma(t_1, t_2)\, = \,\bigl\{\boldsymbol{x}(t) \ | \ 
t \in [t_1, t_2]\bigr\}$ be the arc of the trajectory and  $G(t_1, t_2)\, = \, {\rm co}_s \, \Gamma(t_1, t_2)$, where 
${\rm co}_s X\, =  \, {\rm co} \,\{X , -X\}$ is 
the {\em symmetrized convex hull}. 
It follows easily from the equation~$ \dot{\boldsymbol{x}} = A {\boldsymbol{x}}\, $ that $ \Gamma(t_1+h, t_2+h)\, = \, 
e^{hA}\, \Gamma(t_1, t_2) $ and that for every segment~$ [t_1, t_2] $ the linear span of the arc~$ \Gamma(t_1, t_2) $ coincides with~$L_{{\boldsymbol{x}}_0} $. We also use the short notation   $ \Gamma = \Gamma(0, +\infty) $ for the entire trajectory and, respectively, $ G = G(0, +\infty) $. Thus, $ G $ is the symmetrized convex hull of~$ \Gamma $. 

We denote by $ {\rm ri}\, X $ the relative interior of a set~$ X $, i.e., the interior in its affine span. Every convex subset~$ X\subset \mathbb{R}^d $ that contains more than one point has a nonempty relative interior; if $ X=-X $, then $ 0\in {\rm ri}\, X $. Therefore, $ 0\in {\rm ri}\,  G(t_1, t_2 )$ 
for every segment~$ [t_1, t_2] $. If $L_{\boldsymbol{x}_0} = \mathbb{R}^d$, then $ G(t_1, t_2) $ is a full-dimensional symmetric convex body, and the relative interior becomes the 
usual interior~$ {\rm int}\,  G(t_1, t_2) $.

\begin{defn}\label{d.80}
Let~$ A $ be a Hurwitz $ d\times d $ matrix and $ \dot{ \boldsymbol{x}} \, = \, A\boldsymbol{x}, \ \boldsymbol{x}(0) = \boldsymbol{x}_0$, be the  corresponding linear system. 
A number~$T>0$ is called a {\em cut tail point} for the trajectory~$ \boldsymbol{x}(\, \cdot \,) $
if  for every $ t\ge  T $, the point $ \boldsymbol{x}(t) $  belongs to the relative interior of~$ G $. 
	
A  cut tail point for the matrix~$ A $ is a number~$ T > 0 $ which is a cut tail point of all trajectories.  
\end{defn}   
Since~$ \boldsymbol{x}(t) \to 0 $ as $ t\to +\infty $, and $ 0\in {\rm ri}\, G $, 
it follows that every trajectory of a stable system possesses cut tail points.
The existence of cut tail points of a matrix, i.e., of universal cut tail points of all trajectories,  is less obvious and will be proved below. 
The infimum of all cut tail points of one trajectory~$ \boldsymbol{x}(\, \cdot \,)$ is denoted as~$ T_{cut}(\boldsymbol{x}) $ or as $ T_{cut}(\boldsymbol{x}_0) $. This is the  moment of time when the trajectory enters the interior of its symmetrized convex hull.  
\begin{prop}\label{p.90}
For every  trajectory~$ \boldsymbol{x}(\, \cdot \,) $, the arc $ \Gamma (0, T_{cut}) $, where~$ T_{cut} = T_{cut}(\boldsymbol{x}) $, lies on the (relative) boundary of~$ G$.  
\end{prop}
{\tt Proof.} If, on the contrary, there exists~$ t< T_{cut} $ such 
that~$ \boldsymbol{x}(t) \in {\rm ri} \, G $, it follows that 
~$ \boldsymbol{x}(t) \in {\rm ri} \, G(0,\tau) $ for some~$ \tau>0 $. 
Hence~ $ \boldsymbol{x}(T_{cut})\, = \, e^{(T_{cut}-t) A}\boldsymbol{x}(t)  \in
{\rm ri} \, e^{(T_{cut}-t) A} G(0,\tau)\, = \, {\rm ri} \,  
G(T_{cut}-t, T_{cut}-t + \tau) \subset {\rm ri}\, G $. Thus, 
$ \boldsymbol{x}(T_{cut}) \in  {\rm ri}\, G $, which is a contradiction. 

\begin{figure}
\begin{flushleft}
\includegraphics[height=74mm]{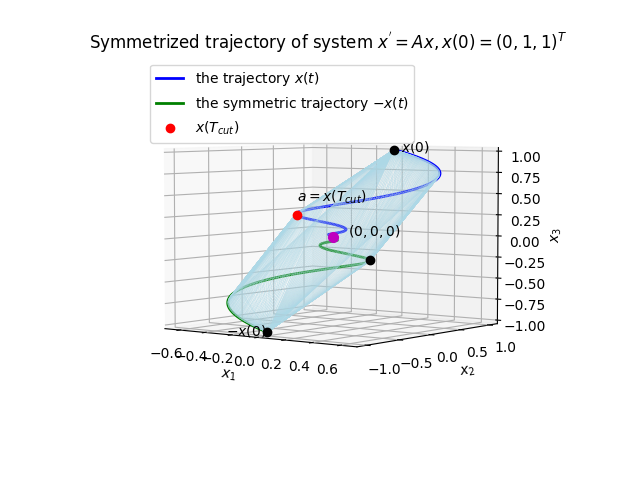}    
\caption{The critical point~$ \boldsymbol{x}(T_{cut})$ on a 3D trajectory}  
\label{fig1}                                 
\end{flushleft}                                 
\end{figure}

We see that each trajectory~$ \boldsymbol{x}(\, \cdot \,) $ consists of two arcs:  
the ``head'' of trajectory $\Gamma(0, T_{cut})$ lies on the  boundary of~$G$ and the tail 
$\Gamma(T_{cut}, +\infty)$  is inside the  interior.
Hence, {\em the set of cut tail points is precisely the half-line
	$(T_{cut}, +\infty)$ and 
	all extreme points of~$G$ are located on the curves~$\, \pm \ \Gamma(0, T_{cut})$}.  
Fig.~\ref{fig1} demonstrates an example in~$\mathbb{R}^3$ with~$\boldsymbol{x}(0) = (0,1,1)$.     

Now let us establish the relation between cut tail points of 
different trajectories. 
\begin{prop}\label{p.100}
	Let $\boldsymbol{x}_0 \in \mathbb{R}^d \setminus \{0\}$; then for every 
	point~$\boldsymbol{y}_0 \in L_{\boldsymbol{x}_0}$, 
	we have~${T_{cut}(\boldsymbol{y}_0) \, \le \, T_{cut}(\boldsymbol{x}_0)}$.  
\end{prop}
{\tt Proof.} Without loss of generality, passing to a subspace, 
it may be assumed that~$L_{\boldsymbol{x}_0} = \mathbb{R}^d$. In this case  the 
vectors~$\{A^j\boldsymbol{x}_0\}_{j=0}^{d-1}$ are linearly 
independent and hence, the vector~$\boldsymbol{y}_0$ can be 
expressed as their linear combination:~$\boldsymbol{y}_0 \, = \, \sum_{j=0}^{d-1}p_jA^j\boldsymbol{x}_0$. 
The operator~$C\, = \, \sum_{j=0}^{d-1}p_jA^j$ commutes with $A$
and $\boldsymbol{y}_0 \, = \, C\boldsymbol{x}_0$. For the function~$\boldsymbol{y}(t)\, = \, C\boldsymbol{x}(t)$, we have  
$\dot{\boldsymbol{y}} \, = \, C \dot{\boldsymbol{x}} \, = \, C A \boldsymbol{x} \, = \, A C \boldsymbol{x}\, = \, A\, \boldsymbol{y}$. 
Thus, $\dot{\boldsymbol{y}} = A\boldsymbol{y}$, and therefore, $\boldsymbol{y}\, = \, C\boldsymbol{x}$ is a trajectory starting at~$\boldsymbol{y}_0$.
The operator $C$ maps the interior of a convex set~$G(0, +\infty)$ to a (relative) interior of the image, 
therefore, if $T$ is a cut tail point for~$\boldsymbol{x}_0$, then so is for~$\boldsymbol{y}_0$.

{\hfill $\Box$}
\medskip 

\begin{cor}\label{c.30}
	If~$L_{\boldsymbol{y}} =  L_{\boldsymbol{x}}$, then  $T_{cut}(\boldsymbol{y}) \, =  \, T_{cut}(\boldsymbol{x})$.  
\end{cor}

Thus, the larger the subspace~$L_{\boldsymbol{x}_0}$ the larger is the value~$T_{cut}(\boldsymbol{x}_0)$.
The biggest possible dimension of~$L_{\boldsymbol{x}_0}$ over all $\boldsymbol{x}_0 \in \mathbb{R}^d$
is  equal to the degree~$n$ of the minimal polynomial of the matrix~$A$. 
We call every subspace~$L_{\boldsymbol{x}_0}$ of dimension~$n$ the {\em largest invariant subspace}
and the corresponding point~$\boldsymbol{x}_0$ {\em generic}. 
If $n< d$, then the largest invariant  subspace is not unique. 
All such subspaces~$L$ have the same dimension~$n$ and the same 
Jordan form and the same Jordan form of the restrictions~$A|_{L}$. It is 
obtained from the Jordan form of~$A$ by leaving only 
one block (the one of the maximal size) of each eigenvalue. 

Applying Proposition~\ref{p.100} and Corollary~\ref{c.30} 
we obtain
\begin{thm}\label{th.40}
	For all generic points~$\boldsymbol{x}$ of a stable matrix~$A$, the 
	value~$T_{cut}(\boldsymbol{x})$ is the same and is the maximal one 
	among all values~$T_{cut}(\boldsymbol{y}), \, \boldsymbol{y}\in \mathbb{R}^d$.  
\end{thm}
As an immediate consequence, we obtain the existence 
of cut tail points of the matrix~$A$: this a cut tail point of  
is an arbitrary generic point~$\boldsymbol{x}$. The common value of~$T_{cut}(\boldsymbol{x})$
for all generic points~$\boldsymbol{x} \in \mathbb{R}^d$ 
will be denoted as~$T_{cut}(A)$. It depends only on 
the minimal polynomial of~$A$.  

\section{Systems with restricted  switching intervals}
Consider a linear switching system with the switching intervals 
restricted not only from below (by the  dwell time) but also from above. 
To every~$A \in \mathcal{A}$ we associate 
a number~$M(A) \, \in \, \bigl(m(A), +\infty \bigr]$. Then $\, \xi = \{\mathcal{A}, \boldsymbol{m}, \mathcal{M} \}$ is 
a linear switching system 
such that   for every regime~$A \in \mathcal{A}$,  the lengths of all switching
intervals corresponding to~$A$ are from the segment~$\bigl[m(A), M(A) \bigr]$.
The systems without upper restrictions (i.e., with~$M=+\infty$) are denoted as 
$\, \xi = \{\mathcal{A}, \boldsymbol{m}\}$.  We use the notation~$\mathcal{T}_{cut}$ for the map
$A \mapsto T_{cut}(A)$. 
\begin{thm}\label{th.10}
	A  linear switching system~$\xi\, =\, \{\mathcal{A}, \boldsymbol{m}\}$ is stable if and only if 
	all its regimes~${A \in \mathcal{A}}$ are stable (i.e., all the matrices~${A}$ are Hurwitz)  and the system~$\xi_{cut}\, =\, \{\mathcal{A}, \boldsymbol{m}, \boldsymbol{m} + \mathcal{T}_{cut}\}$
	is  stable. 
\end{thm}
Thus, if a system with the switching intervals not exceeding~$M = m+T_{cut}$
is stable, then the system with longer switching intervals is stable as well. 
In the proof we use the following result 
from~\cite{PK}. Let to every 
regime 
$A \in \mathcal{A}$, one  associate a norm $f_A$ in~$\mathbb{R}^d$. 
The collection~$f=\{f_A\ | \ A \in \mathcal{A}\}$ is called  a~{\em multinorm}. 
We assume that those norms are uniformly bounded. 
\smallskip 

\smallskip 

\noindent \textbf{Theorem A}~\cite{PK}. {\em A system~$\{\mathcal{A}, \boldsymbol{m}, \mathcal{M}\}$ is stable if and only if 
	there exists a {\em Lyapunov multinorm}~$f$, which possesses the 
	following property: for every trajectory~$\boldsymbol{x}(\cdot)$ and for arbitrary 
	its switching point~$\tau$ between some regimes~$A, A'\in \mathcal{A}$, 
	we have~$f_A(\boldsymbol{x}(\tau))\, > \, \sup\limits_{t \in [\tau + m(A'), \tau + T)]} f_{A'}(\boldsymbol{x}(t))$, 
	where $\tau + T$ is the next switching point}. 
\smallskip

{\tt Proof of Theorem~\ref{th.10}}. The necessity is obvious. 
To prove the sufficiency, we apply~Theorem~A for~$M \, = \, m + T_{cut}$. 
The norm~$f_{A'}(\, \cdot \,)$ is convex, hence, its supremum on the 
arc~$\gamma(\tau + m(A'),\tau + T_{cut}(A') )$ is equal to the supremum on the symmetrized convex 
hull~$G(\tau + m(A'),\tau + T_{cut} )$. 
If we increase $T_{cut}$ to an  arbitrary number~$T > T_{cut}(A)$, 
then the symmetrized convex 
hull does not change, hence neither does the supremum. 
Hence, replacing~$T_{cut}$ by any larger~$T$ does not 
violate the assumptions of Theorem~A and therefore, does not violate the 
stability. Therefore, the system stays stable for~$T=+\infty$.

{\hfill $\Box$}
\smallskip

Thus, if the system with the lengths of all switching intervals 
bounded above by~$\mathcal{T}_{cut}$ is stable, then it stays stable 
without any upper restrictions. In the proof of Theorem~\ref{th.10} we established that 
not only the stability but also the Lyapunov multinorm 
stays the same after removal of all the restrictions. 
\smallskip

\begin{rem}\label{r.20}
	{\em In view of Theorem~\ref{th.10} 
		the restriction of the switching intervals can be realized 
		by finding the critical  point~$T_{cut}(A)$ for every matrix~$A\in \mathcal{A}$. 
		This is, however, not an easy task since it requires  computing the convex hull of the trajectory. 
		To solve this problem it suffices 
		to come up with a method of deciding whether a given~$T>0$ is a cut tail point. 
		Then~$T_{cut}$ can be computed by bisection in~$T$. We show that this decision 
		problem 
		can be reduced to a Chebyshev-like problem of a uniform polynomial approximation.
		That approximation is not by algebraic polynomials, but  by 
		exponential ones (or {\em quasipolynomials}).  To the best of our knowledge, this type of approximation has not been studied in the literature. We present an algorithm of its solution, including the theoretical background. Then we demonstrate the results of numerical implementation. This allows us to compute~$T_{cut}$ for any Hurwitz $d\times d$ matrix and, respectively, to decide the stability of unrestricted systems. 
		
	}
\end{rem}

\section{Characterising~$T_{cut}$ by quasipolynomials of least deviation}
We present a method to decide  whether a given~$T>0$ is a cut tail point
for a matrix~$A$. Then~$T_{cut}(A)$ is computed merely by bisection. 

From the affine similarity of trajectories, it follows that 
for every generic point~$\boldsymbol{x}_0\in \mathbb{R}^d$, the linear space of 
functions~$g_{\boldsymbol{x}_0, \boldsymbol{y}}(t) \, = \, \bigl(e^{tA}\boldsymbol{x}_0, \boldsymbol{y}\bigr),  \  \boldsymbol{y} \in \mathbb{R}^d$, 
is the same. Denote it by~$\mathcal{P}_A$. 
This is the  space of quasipolynomials 
which are linear combinations of functions~$t^ke^{\, \alpha t}
\cos \beta t$ and $t^ke^{\, \alpha t}\, \sin \beta t\, , \ k = 1, \ldots , r-1$, 
where $\lambda = \alpha + i \beta$ is an eigenvalue of~$A$ and  
$r= r(\lambda)$ is the size of the 
largest Jordan block corresponding to~$\lambda$. Thus, the  dimension~${\rm dim}\, {\mathcal{P}}_A\, = \, n$
is equal to the degree of the minimal polynomial of the matrix~$A$. 
\begin{thm}\label{th.70}
	Let~$A$ be a Hurwitz matrix and~$T>0$ be a number; then $T > T_{cut}(A)$ if and only if  
	the value of the minumum in the following  
	extremal problem 
	\begin{equation}\label{eq.extr-cut}
		\left\{
		\begin{array}{l}
			\|\boldsymbol{p}\|_{[0,T]} \ \to \ \min\\
			\boldsymbol{p}(T) = 1 \\
			\boldsymbol{p} \in {\mathcal{P}}_A
		\end{array}
		\right. 
	\end{equation}
	is bigger than~$1$.  
\end{thm}
{\tt Proof}. By Proposition~\ref{p.90}, ~$T >  T_{cut}$ if and only if  
$\boldsymbol{x}(T) \in {\rm int}\, G$. In view of the convex separation theorem, this means that,  for 
every  vector~$\boldsymbol{y} \ne 0$, we have $\bigl(\boldsymbol{y}, \boldsymbol{x}(T)\bigr)\, < \, 
\sup_{\boldsymbol{z} \in G} \bigl(\boldsymbol{y}, \boldsymbol{z}\bigr) \, = \, 
\sup\limits_{t \in [0,T_{cut}]}  \bigl| \bigl(\boldsymbol{y}, \boldsymbol{x}(t)\bigr)\bigr|$. The latter equality 
follows from that~$G = {\rm co}_s\, \Gamma(0, T_{cut})$. 
Since $\boldsymbol{x}(t)\, = \, e^{tA}\boldsymbol{x}_0$, the functions  $\boldsymbol{p}(t)\, = \, \bigl(\boldsymbol{y}, \boldsymbol{x}(t)\bigr)$ 
run over the whole space~$\mathcal{P}_A$. Therefore, $\boldsymbol{p}(T) \, < \, \sup\limits_{t\in [0, T_{cut}]} |\boldsymbol{p}(t)|$
for every~$\boldsymbol{p} \in \mathcal{P}_A$. Replacing~$T_{cut}$ by a bigger~$T$ keeps this inequality. Normalising we get:~$\|\boldsymbol{p}(t)\|_{C[0,T]} > 1$
provided~$\boldsymbol{p}(T)=1$, which concludes the proof.

{\hfill $\Box$}
\medskip 

Problem~(\ref{eq.extr-cut}) means finding the polynomial of least deviation from zero on the segment~$[0,T]$ under the constraint~$\boldsymbol{p}(T) = 1$. We are going to characterize 
this polynomial and find it algorithmically. This will be done 
for a more general problem: for the space~$\mathcal{P}$ of polynomials
$\boldsymbol{p}(t) = \sum_{i=1}^n p_i \varphi_i(t)$ with real coefficients~$p_i$ over an arbitrary 
set of linearly independent continuous real-valued functions~$\{\varphi_k\}_{i=k}^n$ on a compact metric space~$Q$
under the constraint~$\ell(\boldsymbol{p}) = 1$, where~$\ell$ is an arbitrary 
linear functional on~$\mathcal{P}$, find the minimal norm~$\|\boldsymbol{p}\|_{C(Q)}$. In what follows we identify every polynomial~$\boldsymbol{p} \in \mathcal{P}$
with the vector of coefficients~$(p_1, \ldots , p_n)\in \mathbb{R}^n$ and 
a functional~$\ell$ is with the corresponding co-vector in~$\mathbb{R}^n$. 

\begin{rem}\label{r.10}
{\em In many works on linear switching 
systems, the approximation tools have been applied to analyse the trajectories.  It  can be piecewise-linear, polynomial or Pad\'e approximation. 
In our problem, we do not approximate trajectories, but use their precise 
values, which can be expressed by quasipolynomials. 
Then we apply methods  of approximation theory to 
characterize the quasipolynimial of the least deviation from zero (Theorem~\ref{th.70}).}
\end{rem}

\begin{defn}\label{d.20}
	Let $\mathcal{P}$ be the space of polynomials generated by  a set of linearly independent continuous functions 
	$\{\varphi_k\}_{k=1}^n$  on a compact set~$Q$ and  
	$\ell$ be a nontrivial linear functional on~$\mathcal{P}$. The {\em polynomial of least deviation
		from zero} (the {\em Chebyshev-like polynomial})
	corresponding to~$\ell$ is an element~$\boldsymbol{p}\in \mathcal{P}$ of the smallest norm 
	in~$C(Q)$ with the property~$\ell(\boldsymbol{p}) = 1$. 
\end{defn}
The classical Chebyshev polynomials correspond to the case $\varphi_k(t) = t^{k-1}, \, k = 1, \ldots , n$ (algebraic polynomials), 
$Q=[-1, 1]$ and $\ell(\boldsymbol{p})\, = \, 2^{2-n}p_{n}$ (the leading coefficient is~$2^{n-2}$). It is well-known that  for Chebyshev (in other terminology -- Haar) systems~$\{\varphi_k\}_{k=1}^n$ on a segment, when 
each polynomial has at most $n-1$ roots, the Chebyshev-like polynomial exists, unique, 
and possesses similar properties to the standard Chebyshev polynomial.
However, in Theorem~\ref{th.70} we deal with the space of quasipolynomials~$\mathcal{P}_{A}$,
which may not be Chebyshev. In this case the problem is more difficult.

By compactness, for every system $\{\varphi_k\}_{k=1}^n$ and for every functional~$\ell$, 
the Chebyshev-like polynomial exists, although, in general it may not be unique. 
By the convexity
of the norm, all such polynomials have the same norm and form a convex compact subset of~$\mathcal{P}$. For Chebyshev systems,  the Chebyshev alternance theorem holds 
in the classical form. There exists a unique polynomial 
of least deviation from zero, which has an alternance of~$n$ points. 
However, for non-Chebyshev systems, it may not be true, and the problem requires other ideas~\cite{KS}, \cite{SU}.

Thus, the solution of the problem~(\ref{eq.extr-cut}) is a  polynomial of least deviation from zero 
for the space of real-valued exponential polynomials~$\mathcal{P} = \mathcal{P}_{A}$ on the segment $\, Q=[0,T]$ with the functional  $\, \ell(\boldsymbol{p})\, = \, \boldsymbol{p}(T)$. Therefore, deciding whether or not  a given~$T$ is a cut tail point is equivalent  to finding a Chebyshev-like polynomial. 
This will be done in the next section. First of all, we generalize 
the concept of alternance and 
prove an analogue of the Chebyshev - Vall\'e-Poussin theorem~\cite{KS} for characterizing the polynomials 
of least deviation.

We denote the moment vector-function  as $\boldsymbol{u}(t)\, = \, \bigl( \varphi_1(t), \ldots, \varphi_n(t)\bigr) \in \mathbb{R}^n$, and  
for a given polynomial $\boldsymbol{p}\in \mathcal{P}$, let $\, 
\boldsymbol{a}(t)\, =\, [{\rm sign} \, \boldsymbol{p}(t)]\, \boldsymbol{u}(t)$. 
\begin{defn}\label{d.30}
	An alternance of a polynomial 
	$\boldsymbol{p} \in \mathcal{P}$ by a system  $\{\varphi_k\}_{i=k}^n$
	with a functional~$\ell$
	is a collection of points~$\{t_i\}_{i=1}^{N} \subset Q$
	such that the conic hull of vectors~$\boldsymbol{a}(t_i), \, i=1, \ldots , N$
	contains~$\ell$ or~$-\ell$.    
\end{defn}
For a Chebyshev system on a segment, 
Definition~\ref{d.30} gives the classical alternace: $n$ points~$t_1< \cdots < t_{n}$ 
with the alternating signs of~$\boldsymbol{p}(t_i), \, i=1, \ldots , n$.  For non-Chebyshev systems, 
there may be a less number of points and the sign alternating condition 
may fail. It is replaced by 
that the cone spanned by the vectors~$\boldsymbol{a}(t_i)$ contains~$\ell$ or~$-\ell$.   
\begin{thm}\label{th.20}
	A polynomial~$\boldsymbol{p} \in \mathcal{P}$ possesses a least deviation on~$Q$ for a given 
	functional~$\ell$ if and only if it possesses an alternance of 
	at most~$n$ points. 
\end{thm}
{\tt Proof}.  Since the function~$f(\boldsymbol{p}) = \|\boldsymbol{p}\|_{C(Q)}$
is convex on the convex set~$\mathcal{S} = \{\boldsymbol{p} \in \mathcal{P}: \ \ell(\boldsymbol{p})=1\}$
of dimension~$n-1$, it follows that one can invoke the refinement theorem~(see, for instance,~\cite{MIT}), 
and conclude that there exist~$\, s\le n$ points~$\{t_i\}_{i=1}^s$
on~$Q$ such that 
$$
\min_{\boldsymbol{p}\in {\mathcal{S}}}\|\boldsymbol{p}\|_{C(Q)}\ = \ \min_{\boldsymbol{p}\in {\mathcal{S}}}\, \max_{i=1, \ldots , s}|\boldsymbol{p}(t_i)|\, . 
$$ 
The optimal polynomial~$\boldsymbol{p} \in \mathcal{S}$ is characterized by two properties: 1)~$|\boldsymbol{p}(t_i)| = \|\boldsymbol{p}\|, \, i = 1, \ldots , s$
(all points~$t_i$ for which this is not true can be removed); 2)  there is no polynomial~$\boldsymbol{h} \in \mathcal{P}$
such that~$\boldsymbol{p} + \boldsymbol{h}  \in \mathcal{S}$ and  
~$|\boldsymbol{p}(t_i) + \boldsymbol{h}(t_i)| < \|\boldsymbol{p}\|, \, i = 1, \ldots , s$.  
In other terms, the system of strict linear inequalities~$[{\rm sign}\, \boldsymbol{p}(t_i)]\, \boldsymbol{h}(t_i) < 0 $ does not have a solution~$\boldsymbol{h} \in \mathcal{P}$ such that~$\ell(\boldsymbol{h}) = 0$. 
Note that~$[{\rm sign}\, \boldsymbol{p}(t_i)]\, \boldsymbol{h}(t_i)\, = \, \bigl(\boldsymbol{a}(t_i), \boldsymbol{p}\bigr)$. 
Thus, the system of linear inequalities~$\bigl(\boldsymbol{a}(t_i), \boldsymbol{p}\bigr) < 0, \ 
i = 1, \ldots , s$ is incompatible with the 
equation~$\, \bigl(\ell, \boldsymbol{h} \bigr) = 0$. By Farkas' lemma~\cite{MIT}, this is equivalent 
to that~${\rm cone}\, \{\boldsymbol{a} (t_i)\}_{i=1}^n$ (the conic hull) contains either~$\ell$ or~$-\ell$. 

\section{Computing the critical lengths of switching intervals}
We use the 
well-known idea of Remez' algorithm called also the Vall\'e-Poussin procedure, which is  
a  recursive construction of the optimal polynomial, see~\cite{Dz}, \cite{R} for the 
classical case and~\cite{SU} for multivariate generalization. Each iteration  
replaces one point of alternance by a point of the maximal deviation of the current polynomial
from zero followed by an averaging of all deviations over the new alternance. 
Since we work with a polynomial under a linear constraint, we have to 
modify this principle. This is done in the  following 

\noindent \textbf{Algorithm~1}. {\em Initialisation.} 
Take an arbitrary set~$\mathcal{T}_1\ = \, \{t_i\}_{i=1}^n$ of $n$ different points on~$Q$ 
and consider the set of moment vectors~$\{\boldsymbol{u}(t_i)\}_{t_i\in \mathcal{T}_1}$. 
We assume they are linearly independent.  
Choose the signs~$\sigma_i \in \{1, -1\}$ so that 
the conic hull~$\mathcal{K}_1$ of~$n$ vectors~$\boldsymbol{a}_i = \boldsymbol{a}(t_i) = \sigma(t_i)\boldsymbol{u}(t_i), \, t_i \in \mathcal{T}_1$, contains the vector~$\ell$, i.e., $\ell = 
\sum_{i=1}^n \alpha_i \boldsymbol{a}_i$ with positive~$\alpha_i$. 
Find the solution~$\boldsymbol{p}_1 \in \mathbb{R}^n, \, r_1 \in \mathbb{R}$
of the linear system~$(\boldsymbol{p}_1, \ell) = 1, \, (\boldsymbol{p}_1, \boldsymbol{a}_{i}) = b_1, 
\ i =1, \ldots , n$.  Since $\ell \in \mathcal{K}_1$, we have $b_1 > 0$. 

Set~$B_1 = \|\boldsymbol{p}_1\|$.  Choose arbitrary~$\varepsilon > 0$. 
\medskip 

{\em Main loop.} {\em The $k$th iteration}. 
We have a set of~$n$ points~$\mathcal{T}_{k}\, = \, \{t_i\}_{i=1}^n$, a polynomial~${\boldsymbol{p}_k \in \mathcal{P}}$, 
and two numbers~$B_k > b_k> 0, \, \sigma \in \{-1, 1\}$
such that: 1) $\ \boldsymbol{p}_k(t_i)\, = \, \sigma_i b_k$, where $\sigma_i \in \{-1, 1\}$ for all~$i$; 
2)  the cone~$\mathcal{K}_k \, = \, {\rm cone}\, \{\boldsymbol{a}_i\}_{i=1}^n$ 
contains the vector~$\sigma\, \ell$, where 
$\boldsymbol{a}_i\, = \, \boldsymbol{a}(t_i)\, = \, \sigma_{i}\boldsymbol{u}(t_i)$.

Denote~$r_k = ~\|\boldsymbol{p}_{k}(\,\cdot \,)\|_{C(Q)}$
and find a point~$t_0$ for which~$|\boldsymbol{p}_{k}(t_0)| \, = \, r_{k}$. 
Clearly, $r_k \ge b_k$. Denote by~$\sigma_0$ the sign of~$\boldsymbol{p}_{k}(t_0)$ and
$\boldsymbol{a}_0 = \boldsymbol{a}(t_0) = \sigma_0\boldsymbol{u}(t_0)$. Denote by~$(\alpha_i)_{i=1}^n$
and~~$(s_i)_{i=1}^n$ the coordinates of the vectors~$\sigma \, \ell$ and $\boldsymbol{a}_0$
respectively in the basis~$\{\boldsymbol{a}_i\}_{i=1}^{n}$. All the numbers~$\alpha_i$
are positive.  

Let the maximal ratio~$\bigl\{ \frac{s_i}{\alpha_i}: \, i=1, \ldots , n\bigr\}$ be attained at some index~$i=q$. 
Define~$\mathcal{T}_{k+1}$ as the set~$\mathcal{T}_{k}$ with the element~$t_q$ replaced by~$t_0$ and 
set~$\gamma_0 = \bigl|\frac{\alpha_q}{s_q}\bigr|$ and $\gamma_j \, = \, \bigl|\alpha_j - \frac{\alpha_j s_j}{s_q}\bigr|$
for~$j\notin \{0,q\}$. 
Then we have ${\rm sign}\, (s_q)\, \sigma \, \ell \, = \, \sum_{j\ne q} \gamma_i \boldsymbol{a}_i$. 

Find the solution~$\boldsymbol{p} \in \mathbb{R}^n, \, r  \in \mathbb{R}$
of the linear system~$(\boldsymbol{p} , \ell) = 1, \, (\boldsymbol{p} , \boldsymbol{a}_{i}) = r, 
\ i =0, \ldots , d, \ i\ne q$, and denote it as~$\boldsymbol{p}_{k+1} = \boldsymbol{p}, \,b_{k+1} = r$. 
Set~$\ B_{k+1}\, = \, \min\, \{r_{k}, B_k\}$. 

Finally, redenote the points of the set~$\mathcal{T}_{k+1}$ by $\{t_i\}_{i=1}^d$
in an arbitrary order.

{\em Termination.} If~$B_k - b_k < \varepsilon$, then STOP. 
The value of the problem~(\ref{eq.extr-cut}) and the norm of the 
Chebyshev-like polynomial~$\boldsymbol{p}$ belongs to the segment~$[b_k, B_k]$.  
The approximation of the polynomial of the least deviation is~$\boldsymbol{p}_k$. 
\smallskip

\begin{rem}\label{r.25}
	{\em  For solving problem~(\ref{eq.extr-cut}) we 
		apply Algorithm~1 to the quasipolynomials~$\mathcal{P} = \mathcal{P}_A$, on the set
		$Q = [0, T]$ and $\ell(\boldsymbol{p}) = \boldsymbol{p}(T)$, i.e., for $\ell = \boldsymbol{u}(T)$. 
		If after $k$th iteration we obtain~$b_k > 1$, then $T$ is a cut tail point, i.e., 
		$T > T_{cut}$; if~$\, B_k \le 1$, then~$T \le T_{cut}$. 
	}
\end{rem}

Now we estimate the rate of convergence of Algorithm~1. 
Denote ${\Gamma = \sum\limits_{i\ne q}\gamma_i}$. 
\begin{thm}\label{th.30}
	In the $k$th iteration of Algorithm~1, we have 
	\begin{equation}\label{eq.estim}
		B_{k+1}-b_{k+1} \ \le  \ \Bigl( 1 - \frac{\gamma_0}{\Gamma}\Bigr)\, 
		\Bigl(B_{k} \ - \ b_{k}\Bigr). 
	\end{equation}
\end{thm}
{\tt Proof}. We have $(\boldsymbol{p}_k, \boldsymbol{a}_i)\, = \, |\boldsymbol{p}_k(t_i)|\, = \, b_k$
for all~$i=1, \ldots , n\, , \ i\ne q$ and 
$(\boldsymbol{p}_k, \boldsymbol{a}_0)\, = \, |\boldsymbol{p}_k(t_0)|\, = \, r_k$. 
On the other hand~$\ (\boldsymbol{p}_{k+1}, \boldsymbol{a}_i)\, = \, |\boldsymbol{p}_{k+1}(t_i)|\, = \, b_{k+1}, \ 
i\ne q$. 
Taking the difference, we obtain~$(\boldsymbol{p}_{k+1} - \boldsymbol{p}_k, \boldsymbol{a}_i)\, = \, 
b_{k+1}-b_k$ for all~$i=1, \ldots , n\, , \ i\ne q$ and 
$(\boldsymbol{p}_{k+1} - \boldsymbol{p}_k, \boldsymbol{a}_0)\, = \,  b_{k+1}-r_k$. 
Multiplying this equality by~$\gamma_i$, taking the sum over all~$i\ne q$, 
we get 
$$
\sum_{j\ne q} \, \gamma_i \, \bigl(\boldsymbol{p}_{k+1} - \boldsymbol{p}_k\, , \, \boldsymbol{a}_i\bigr) \ =\ 
\bigl(\boldsymbol{p}_{k+1} - \boldsymbol{p}_k\, , \, {\rm sign}\, (s_q)\, \sigma \, \ell \bigr) \ = \
$$
$$ 
= \ {\rm sign}\, (s_q)\, \cdot \, \sigma \, \cdot \, \Bigl[ \bigl(\boldsymbol{p}_{k+1}\, , \, \ell \bigr)\ - \ 
\bigl(\boldsymbol{p}_k\, ,  \, \ell \bigr)\Bigr]
\ = \ 0\, , 
$$

because $\ (\boldsymbol{p}_{k+1}\, ,  \, \ell )\, = \, (\boldsymbol{p}_k\, ,  \, \ell )\, = \, 1$. 
Therefore, 
$$
0 \, = \,  
\sum_{j\ne q} \, \gamma_i \, \bigl(\boldsymbol{p}_{k+1} - \boldsymbol{p}_k\, , \, \boldsymbol{a}_i\bigr) \, =\,   
\gamma_0\, \bigl( b_{k+1}-r_k\bigr) \,  + \,  
$$
$$
+ \, \sum_{j \notin  \{0, q\}} \, \gamma_i \, \bigl(b_{k+1}-b_k\bigr) \, = \, \gamma_0\,  \bigl( b_{k}-r_k\bigr)  + \sum_{j \ne q} \, \gamma_i \, \bigl(b_{k+1}-b_k\bigr).  
$$
Therefore, 
$\gamma_0\, ( r_k - b_{k}) \ = \  \Gamma \, (b_{k+1}-b_k)$, and hence, 
$$
b_{k+1}-b_k \ = \ \frac{\gamma_0}{\Gamma}\, \Bigl(r_k \ - \ b_{k}\bigr)\ \ge  \ 
\frac{\gamma_0}{\Gamma}\, \Bigl(r_{k} \ - \ b_{k}\Bigr). 
$$
Rewriting we obtain
$$
r_{k}-b_{k+1} \ \le  \ \Bigl( 1 - \frac{\gamma_0}{\Gamma}\Bigr)\, 
\Bigl(r_{k} \ - \ b_{k}\Bigr). 
$$
Let us remember that~$B_{k+1} \, = \, \min\, \bigl\{ B_k, r_k \bigr\}$. 
If $r_k \le  B_k$, then we replace in the left hand side $r_k = B_{k+1}$, 
increase the right hand side by substituting~$B_k$ for $r_k$  and obtain~(\ref{eq.estim}). 
If~$r_k >  B_k$, then we keep the inequality true 
by reducing $r_k$ to $B_k$, then replace $B_k = B_{k+1}$ in the left hand side and 
arrive at~(\ref{eq.estim}).

{\hfill $\Box$}
\medskip 

\smallskip 

\begin{cor}\label{c.20}
	If in all iterations of Algorithm~1, the parameter~$\frac{\gamma_0}{\Gamma}$
	is not less than a number~$\mu \, >\, 0$, then the algorithm has a linear convergence and
	$$
	B_{k+1} \, - \, b_{k+1} \ \le \ (1-\mu)^k \, \bigl(B_{1} - b_{1}\bigr)\, , \quad k \in \mathbb{N}\, . 
	$$
\end{cor}

\begin{rem}\label{r.35}
	{\em By Corollary~\ref{c.20}, Algorithm~1 converges linearly provided 
		the numbers~$\gamma_0$ in all iterations are uniformly bounded below by some constant~$\mu> 0$. 
		By definition, $\gamma_0$ is proportional to the coefficient~$\alpha_q$
		in the expansion~$\boldsymbol{a}_0 \, = \, \sum_{i=1}^n \alpha_i\boldsymbol{a}_i$. If~$\alpha_q$ 
		is small,  then the convergence  may slow down. 
		This happens when the vector~$\boldsymbol{a}_0$ is close 
		to the hyperplane spanned by the~$d-1$ vectors~$\boldsymbol{a}_i, \, i\ne q$. 
		Equivalently, the $n$ vectors~$\boldsymbol{u}(t_0), \ldots , \boldsymbol{u}(t_n)$ (without~$\boldsymbol{u}(t_q)$)
		are almost linearly dependent. In practice, if this singularity phenomenon occurs at some iteration, 
		one should slightly perturb the point~$t_0$, without significant descent 
		of the value~$|\boldsymbol{p}(t_0)|$.  
		
		If all exponents of the quasipolynomials are real, i.e., 
		if the spectrum of the matrix~$A$ is real, 
		then these exponents form a Chebyshev system~\cite{KS} and the 
		problem~(\ref{eq.extr-cut}) can be solved by the Remez algorithm~\cite{R}.
		The singularity phenomenon is impossible in this case, and the Remez algorithm
		always converges linearly~\cite{Dz}.} 
\end{rem}

\section{Numerical results and examples}
We apply Algorithm~1 for computing~$T_{cut}(A)$ and deciding the stability 
of linear switching systems. We begin with a two-dimensional example.  

\noindent \textbf{Example~1}. 
Consider a linear switching system $ \xi = (\mathcal{A}, m, M) $, with $ m = 1, M = 2.5 $ 
and with the control set~$\mathcal{A} = \{A_1, A_2\}$, where

\small 
\[ 
A_1 = \left( 
\begin{array}{cc}
	-0.3216 & -1 \\
	2 & -0.3216
\end{array} \right), \,
A_2 = \left( 
\begin{array}{cc}
	-0.3216 & -2 \\
	1 & -0.3216
\end{array} \right).
\]

Since~$M=2.5m$, we see that the system has pretty often switches.  	
Applying the invariant polytope method from~\cite{PK} we 
get that the Lyapunov exponent is negative and bigger than~$-0.003$. 
Thus, the system is stable but is very close to instability since the 
Lyapunov exponent is close to zero. What will happen if we remove the upper restriction
$M$, i.e., allow arbitrary switching intervals longer than~$m=1$? Will the system stay stable? Intuitively it should loose the stability, but in reality it does not. 
Applying Algorithm $ 1 $ we have $ T_{cut}(A_1) = T_{cut}(A_2) = 1.4239 $. Since 
$ m + T_{cut}  = 1.4239 <  M$,  it follows from Theorem~\ref{th.20} 
that the system~$\{\mathcal{A}, m \}$ is stable.

\noindent  \textbf{Computation of $ T_{cut}(A) $ in different dimensions}.

We apply Algorithm~1 for computing~$T_{cut}(A)$
and the corresponding critical time intervals 
in dimensions~$d$  from~2 to 12.  The results are given in Table~\ref{tab1}. 

\begin{table}[h!]
	\centering
	\caption{$ T_{cut}(A) $ in dimensions~$d$  from~2 to 12.}
	\label{tab1}
	\begin{tabular}{|c|c|c|}
		\hline
		$  d $ & $ T_{cut}(A) $ & $ t $ 
		\\ 
		\hline
		$ 2 $ & $ 9.0133 $  & $ 17 $
		\\ 
		\hline
		$ 4 $ & $ 13.7237 $ & $ 52 $
		\\ 
		\hline
		$ 6  $ & $ 72.4608 $  & $ 149 $
		\\ 
		\hline
		$ 8 $ & $ 85.6227 $ & $ 233 $
		\\ 
		\hline	
		$ 10 $ & $ 177.1417 $ & $ 408 $
		\\ 
		\hline
		$ 12 $ & $ 47.0105 $ & $ 969 $
		\\ 
		\hline
	\end{tabular}
\end{table}

We run Algorithm~$ 1 $ until $ B_k - b_k \leqslant 10^{-4} $, so~$ T_{cut}(A) $
is found with the accuracy~$ 10^{-4}$. The third column shows the computational time~$t$. 
We see that it takes a few seconds for the dimensions~$d\le 4$
and about a quarter  of an hour for~$d=12$. 

\section{Planar systems}
For $2\times 2$ matrices, the parameter $T_{cut}$ can be evaluated geometrically and in an  explicit form. 
We have a real $2\times 2$-matrix~$A$. 
For the sake of simplicity we exclude the case of multiple eigenvalues.

\begin{figure}
	\begin{center}
		\includegraphics[height=3.9cm]{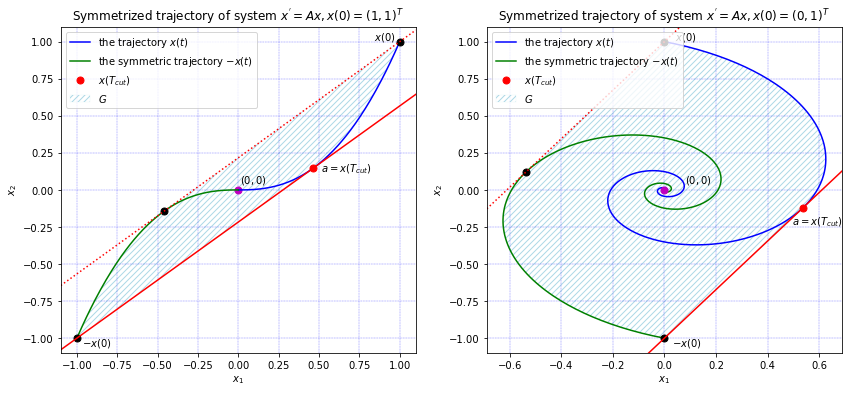}    
		\caption{The point~$\boldsymbol{x}(T_{cut})$ on a 2D trajectory. Left: the case of real 
			eigenvalues; Right: complex eigenvalues}  
		\label{fig2}                                 
	\end{center}                                 
 \end{figure}

\noindent {\tt Case 1. $A$ has real $\alpha_1 < \alpha_2 <0$}. 
Denote them by $\alpha_1, \alpha_2$. 
The trajectory~$\Gamma =  \{\boldsymbol{x}(t)\, | \, t \ge 0\}$ in the basis of eigenvectors 
is~$(x_1(t),  x_2(t)) \, = \, 
\bigl( e^{\alpha_1 t} \, , \,   e^{\alpha_2 t}\bigr)\, , \ 
t \in \mathbb{R}_+$. 
The tangent line 
from the point~$(-1, -1)$ touches~$\Gamma$ at the point~$\boldsymbol{a} = \boldsymbol{x}(T_{cut})\in \Gamma$
(fig.~\ref{fig2} left), 
where $t=T_{cut}$ is found from the system of equations
(we denote~$\, \boldsymbol{a} \, + \, s \, \dot{\boldsymbol{x}}(T_{cut}) \, = 
\, (-1,-1)$):
\begin{equation}
	\left\{
	\begin{array}{lcl}
		e^{\alpha_1 t}\ + \ s\, \alpha_1\, e^{\alpha_1 t}& = & - 1\\
		e^{\alpha_2 t}\ + \ s\, \alpha_2\, e^{\alpha_2 t}& = & - 1
	\end{array}
	\right. 
\end{equation}
which becomes after simplification 
$\frac{1 + e^{- \alpha_1 t}}{\alpha_1}\, = \, \frac{1 + e^{- \alpha_2 t}}{\alpha_2}$. 
\medskip

{\tt Case~2. Two complex eigenvalues}
~$\alpha \pm i \beta$ with $\alpha < 0$. 
The trajectory~$\Gamma$ in a suitable basis 
has the equation~$(x_1(t), x_2(t)) \,  = \, 
e^{\alpha t}\, (\cos \beta t \, , \, \sin \beta t)
\, , \ 
t \in \mathbb{R}_+$. The trajectory $\Gamma$
goes from the point~$\boldsymbol{x}(0) = (1, 0)$ to zero making infinitely  many rotations. 
Taking the point of tangency~$\boldsymbol{a}$ of $\Gamma$ with the line 
going from the point~$(-1, 0)$, fig.~\ref{fig2} right,  we obtain 
\begin{equation}\label{eq.prec2}
	\left\{
	\begin{array}{lcr}
		e^{\alpha t} \cos \beta t \ + \ s\, e^{\alpha t}\bigl(\alpha \cos \beta t \, - \, 
		\beta \sin \beta t \bigr)& = & - 1\\
		e^{\alpha t} \sin  \beta t \ + \ s\, e^{\alpha t}\bigl(\alpha \sin \beta t \, + \, 
		\beta \cos \beta t \bigr)& = & 0
	\end{array}
	\right. 
\end{equation}
from which it follows  $\alpha \sin \beta t \, + \, \beta \cos \beta t\, + 
\, \beta \, e^{\alpha t} \, = \, 0$. The unique solution of this equation is $t = T_{cut}$. 
\vspace{1cm} 

\noindent \textbf{Acknowledgements}. 
The work  is performed with the support of the Theoretical Physics and Mathematics Advancement Foundation ``BASIS''.

\end{document}